# Effects of Improved-Floor Function on the Accuracy of Bilinear Interpolation Algorithm

Olivier Rukundo

Department of Communication and Information Sciences, Tilburg University, Tilburg, Netherlands

Correspondence: Olivier Rukundo, Department of Communication and Information Sciences, Tilburg University, Tilburg, PO Box 90153 5000 LE Tilburg, Netherlands. Tel: +31 13 466 8118. E-mail: orukundo@gmail.com

**Abstract**

In this study, the standard IEEE 754–2008 and modulo-based floor functions for rounding non-integers have been presented. Their effects on the accuracy of the bilinear interpolation algorithm have been demonstrated. The improved-floor uses the modulo operator in an effort to make each non-integer addend an integer when the remainder between the multiplicative inverse of the fractional factor and the numerator is greater than zero. The experiments demonstrated relatively positive effects with the improved-floor function and the alternativeness to the standard round function.

## 1. Introduction

In many computational fields, the bid for optimal accuracy remains an unsolved problem. For example, improving the accuracy of spatially image interpolation algorithms generally involves the use of many adjacent pixels to compute the missing pixel values in high resolution (H.R) images. However, the processing time involved is more when adjacent pixels are used. Furthermore, the low-computational complexity emerges as another important requirement, since in many cases, the wide-spread use depends on how accurately fast and low-cost the algorithms are (Rukundo & Cao, 2012, pp. 25-30, Ramponi, 1999, pp. 629-639, Rukundo & Maharaj, 2014, pp. 641-647). In fact, this arbitrary selection or use of just 'many adjacent pixels' does not even consider the negative effects of the numerical errors on the overall performance of the numerical scheme. Therefore, this study focuses on the effects of the numerical errors, particularly the effects of the improved-floor function against the standard round-off functions on the bilinear image interpolation algorithm. The bilinear algorithm has an advantage of producing reasonably realistic images with only four adjacent pixels, at an acceptable computational complexity compared to other linear function-based image interpolation algorithms. In this study, the finite precisions of computations are not the main concern. However, only the negative effects of the round-off errors have been studied (Higham, 1996). Therefore, the improved floor rounding scheme or function has been developed to adapt to the cases with relatively large or small number of half-integers in H.R image interpolation. Note that, a half-integer is a number of the form $n + 1/2$, where $n$ is an integer. Previously, the standard IEEE 754–2008 defined five rounding rules in the homogeneous general computational operations, as shown with examples of only half-integers in Table 1.

Table 1. Five rounding rules defined by the IEEE 754-2008 standard.

| MODES | | EXAMPLES OF HALF-INTEGERS | | | |
|---|---|---|---|---|---|
| | | +11.5 | +12.5 | -11.5 | -12.5 |
| round to nearest, ties/half-integers to even | | +12.0 | +12.0 | -12.0 | -12.0 |
| round to nearest, ties/half-integers away from zero | (*round*) | +12.0 | +13.0 | -12.0 | -13.0 |
| Round toward 0 | (*fix*) | +11.0 | +12.0 | -11.0 | -12.0 |
| round toward +∞ | (*ceil*) | +12.0 | +13.0 | -11.0 | -12.0 |
| round toward −∞ | (*floor*) | +11.0 | +12.0 | -12.0 | -13.0 |

Subsequently, Maxfield published, in his book and website, a simplified diagram of main rounding modes (see Table 2) and he also included the examples of 'non-half-integers' (Maxfield & Brown, 2005). Here, a non-half-integer is a non-integer other than a half-integer. According to Maxfield, the round-toward-nearest encompasses both the R-H-U and R-H-D modes. For arithmetic Rounding or R-H-U, there can be both symmetric and asymmetric versions. Round-up may be equated to R-C or R-A-T-Z depending on the application/implementation. Truncation is identical to R-T-Z. However, in any case, the rule or mode perceived the best and has always been dependent on the number of half-integers encountered in the rounding process. Therefore, each rule has been working better in one case and worse in the other; thus raising the problematic question about a scheme that would best fit in all the cases (yielding results perceivable optimal). In this effort, the modulo operator has been used to deal with the whole and fractional factor of each addend, in the conventional bilinear interpolation equation (1). Furthermore, the obtained sum is technically rounded towards negative infinity. Since the dimensions between the reference and interpolated images must be identical at various scaling ratios, in the experiments, the reference images were downsampled to match the dimensions of the interpolated images, and the ACDsee software has been used for the respective scenarios.

Table 2. Maxfield's diagram

| Mode | ← Negative Infinity | | | | | | | Zero → | ← Zero | | | | | Positive Infinity → | | |
|---|---|---|---|---|---|---|---|---|---|---|---|---|---|---|---|---|
|  | -2.0 | -1.7 | -1.5 | -1.3 | -1.0 | -0.7 | -0.5 | -0.3 | 0.0 | 0.3 | 0.5 | 0.7 | 1.0 | 1.3 | 1.5 | 1.7 | 2.0 |
| R-H-U (s) | -2 | -2 | -2 | -1 | -1 | -1 | -1 | -0 | 0 | 0 | 1 | 1 | 1 | 1 | 2 | 2 | 2 |
| R-H-U (a) | -2 | -2 | -1 | -1 | -1 | -1 | -0 | -0 | 0 | 0 | 1 | 1 | 1 | 1 | 2 | 2 | 2 |
| R-H-D (s) | -2 | -2 | -1 | -1 | -1 | -1 | -0 | -0 | 0 | 0 | 0 | 1 | 1 | 1 | 1 | 2 | 2 |
| R-H-D (a) | -2 | -2 | -2 | -1 | -1 | -1 | -1 | -0 | 0 | 0 | 0 | 1 | 1 | 1 | 1 | 2 | 2 |
| R-H-E | -2 | -2 | -2 | -1 | -1 | -1 | -0 | -0 | 0 | 0 | 0 | 1 | 1 | 1 | 2 | 2 | 2 |
| R-H-O | -2 | -2 | -1 | -1 | -1 | -1 | -1 | -0 | 0 | 0 | 1 | 1 | 1 | 1 | 1 | 2 | 2 |
| R-C | -2 | -1 | -1 | -1 | -1 | -0 | -0 | -0 | 0 | 1 | 1 | 1 | 1 | 2 | 2 | 2 | 2 |
| R-F | -2 | -2 | -2 | -2 | -1 | -1 | -1 | -1 | 0 | 0 | 0 | 0 | 1 | 1 | 1 | 1 | 2 |
| R-T-Z | -2 | -1 | -1 | -1 | -1 | -0 | -0 | -0 | 0 | 0 | 0 | 0 | 1 | 1 | 1 | 1 | 2 |
| R-A-F-Z | -2 | -2 | -2 | -2 | -1 | -1 | -1 | -1 | 0 | 1 | 1 | 1 | 1 | 2 | 2 | 2 | 2 |

*Note*. R-H-U = Round-Half-Up; R-H-E = Round-Half-Even, R-C = Round Ceiling, R-T-Z = Round-Toward-Zero, R-H-D = Round-Half-Down, R-H-O = Round-Half-Odd, R-F = Round Floor, R-A-F-Z = Round-Away-From-Zero, (s) = Symmetric, (a) = Asymmetric.

One of the best structural and pixel based IQA metrics as well as the linear correlation metrics have been used to quantify the effects of the compared rounding schemes as well as the Matlab *tic-toc* commands for run-time overhead (Wang, Bovik & Lu, 2002, pp. IV-3313 – IV-3316, Martin & Bill, 2011). This paper is organized as follows. Section 1 describes the background followed by section 2 that demonstrates the non-integerness found in the process of bilinear interpolation. Section 3 presents the improved floor rounding scheme. Section 4 gives the experimental results and section 5 gives the conclusion.

**2. Non-integer cases in bilinear interpolation**

Numerous studies that clearly explain how the bilinear interpolation algorithm works and that demonstrate its wide application in the literature exist. Therefore, the redundant explanations and demonstrations are not included in this study. The mathematical expression for the bilinear interpolation is shown in equation (1). This equation is used to calculate the H.R. pixel values that have always been either of integer or non-integer type.

$$J(r',c') = I(r,c) \times (1-\Delta r) \times (1-\Delta c) + I(r+1,c) \times \Delta r \times (1-\Delta c) \\ + I(r,c+1) \times (1-\Delta r) \times \Delta c + I(r+1,c+1) \times \Delta r \times \Delta c \quad (1)$$

where, $r$, $c$, $\Delta r$ and $\Delta c$ represent the distances between the pixels, as shown in Figure 1; $J(r',c')$ represents the H.R. pixel values (that need round-off, if non-integers); $I(r,c)$, $I(r+1,c)$, $I(r,c+1)$ and $I(r+1,c+1)$ represent the low resolution (L.R) pixel values; and $(1-\Delta r) \times (1-\Delta c)$, $\Delta r \times (1-\Delta c)$, $(1-\Delta r) \times \Delta c$ and $\Delta r \times \Delta c$ are linearly distributed weights. In this study, the improvement made was based on turning the equation (1)'s non-integer-addends into integers. Figure 1 illustrates how the pixel values are calculated at the H.R. locations shown in red, green, purple, yellow and orange. Further details can be found in (Rukundo & Maharaj, 2014, pp. 641-647). As mentioned earlier, there is a possibility of the periodic appearances of the non-integer cases according to the type of data being processed. Tables 3, 4, 5 and 6 show the circumstances in which such cases of half-integers, non-half-integers and integers occur. The existing strategies do not effectively adapt to such periodic appearances in a manner that is unquestionable. Furthermore, this encourages the need for finding the most effective round-off operation capable of effectively reducing the negative effects of such cases on the accuracy of bilinear image interpolation algorithm. Note that, in this manner, the final objective remains making the interpolated image quality as much artefact-free as possible. The results shown in Tables 3, 4, 5 and 6 have been obtained using the equation (1) with reference to the coordinate values mentioned.

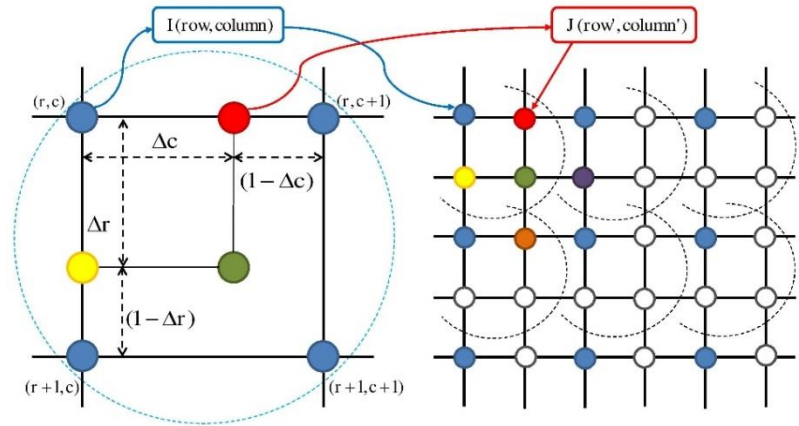

Figure 1. I (row, column) describe generally the L.R image pixels. (r, c), (r+1, c), (r, c+1) and (r+1, c+1) represent the L.R pixels, in blue circles. J (row', column') generally represents the interpolated H.R pixels. In this figure, they are shown in red, yellow and green circles. Their values are systematically calculated using the four L.R. pixels, encircled by a dashed-blue circle

Table 3. With $I(r,c)=91$, $I(r+1,c)=162$, $I(r,c+1)=210$, $I(r+1,c+1)=95$

| $\Delta r$ | $\Delta c$ | $J(r',c')$ |
|---|---|---|
| 0 | 0.5 | 150.5000 |
| 1 | 0.5 | 128.5000 |
| 0.5 | 1 | 152.5000 |
| 0.5 | 0 | 126.5000 |
| 0.5 | 0.5 | 139.5000 |

| 91 | 150.5 | 210 |
|---|---|---|
| 126.5 | 139.5 | 152.5 |
| 162 | 128.5 | 95 |

Table 3 shows the four arbitrary-selected pixels/values to represent the L.R. image pixels. Their locations can be approximated as shown in Figure 1. The values shown in the squares coloured in red, yellow, orange and purple are interpolated values at those H.R. locations. These H.R. values are all half-integers. This implies that the effects of the round-off errors would likely be higher than in any other cases.

Table 4. With $I(r,c)=125$, $I(r+1,c)=99$, $I(r,c+1)=255$, $I(r+1,c+1)=17$

| $\Delta r$ | $\Delta c$ | $J(r',c')$ |
|---|---|---|
| 0 | 0.5 | 190.0000 |
| 1 | 0.5 | 58.0000 |
| 0.5 | 1 | 136.0000 |
| 0.5 | 0 | 112.0000 |
| 0.5 | 0.5 | 124.0000 |

| 125 | 190 | 255 |
|---|---|---|
| 112 | 124 | 136 |
| 99 | 58 | 17 |

In Table 4, the four arbitrary selected pixel values are different from those in Table 3. The selected values have led to integer output at all the H.R pixel locations. This implies that the negative effects of round-off errors will not exacerbate the calculation-errors in the same manner as in the different cases. Furthermore, this situation would to some extent save the interpolated data. In Table 5, the four arbitrary-selected pixels/values are different from those in Tables 3 and 4. The output obtained comprises half-integers, integers and non-half-integers-with-the-fractional-part greater than a half. The round-off errors would not be at the same degree as if only half-integers were obtained. In Table 6, the four arbitrary selected pixels/values differ from those in the previous tables in type by a non-half-integer-with-the-fractional-part less than a half. The round-off effects will not be the same depending on the rounding scheme used. Note that a half of an integer is not always a half-integer. For example, a half of an even integer is an integer but not a half-integer (Sabin, 2010, pp. 218).

Table 5. With $I(r,c) = 191$, $I(r+1,c) = 102$, $I(r,c+1) = 111$, $I(r+1,c+1) = 195$

| $\Delta r$ | $\Delta c$ | $J(r',c')$ |
|---|---|---|
| 0 | 0.5 | 151.0000 |
| 1 | 0.5 | 148.5000 |
| 0.5 | 1 | 153.0000 |
| 0.5 | 0 | 146.5000 |
| 0.5 | 0.5 | 149.7500 |

Table 6. With $I(r,c) = 32$, $I(r+1,c) = 33$, $I(r,c+1) = 72$, $I(r+1,c+1) = 72$

| $\Delta r$ | $\Delta c$ | $J(r',c')$ |
|---|---|---|
| 0 | 0.5 | 52.0000 |
| 1 | 0.5 | 52.5000 |
| 0.5 | 1 | 72.0000 |
| 0.5 | 0 | 32.5000 |
| 0.5 | 0.5 | 52.2500 |

## 3. The Floor Improvement Scheme

As mentioned earlier, the modulo operator is applicable to the fractional factors and to apply the modulo operator it is necessary to simplify equation (1) to clearly distinguish between the pixel (integer) and weight (integer or non-integer) parts, as shown in the equation (2).

$$J(r',c') = N1 \times W1 + N2 \times W2 + N3 \times W3 + N4 \times W4 \qquad (2)$$

where $I(r,c) = N1$, $I(r+1,c) = N2$, $I(r,c+1) = N3$, $I(r+1,c+1) = N4$ and $(1-\Delta r) \times (1-\Delta c) = W1$, $\Delta r \times (1-\Delta c) = W2$, $(1-\Delta r) \times \Delta c = W3$ and $\Delta r \times \Delta c = W4$. With equation (2), the integer factor ($N$) is clearly distinguishable from the weight factor ($W$), and therefore can be dealt with separately. Given that the weight factor ($W$) can be more non-integer (many times) than of integer type, one can make it a multiplicative inverse to maximize chances of using the denominator of integer type in the calculations. Furthermore, with reference to the $\Delta r$ and $\Delta c$ values shown in Tables 3, 4, 5 and 6, one of the two integers $W$ can be zero, and if $W = 0$, this would lead to the 'undesired' indeterminate case; therefore, a small positive number $L$ is added to the denominator of multiplicative inverse of the weight factor ($W$).

$$J(r',c') = \frac{N1}{\left(\frac{1}{W1}\right)} + \frac{N2}{\left(\frac{1}{W2}\right)} + \frac{N3}{\left(\frac{1}{W3}\right)} + \frac{N4}{\left(\frac{1}{W4}\right)} \qquad (3)$$

With reference to equation (3), the first addend can be written as $\frac{N1}{A}$. Now, considering the small positive number $L$, the equivalent of $A$ can be written as $A = \frac{1}{W1+L}$. Applying the modulo operator between the dividend $N1$ and divisor $A$ does effectively determine whether the product will be integer or non-integer depending on the remainder obtained. For example, if $(N1)\mod(A) = 0$ then following addend

$\dfrac{N1+((N1)\bmod(A))}{A}$ would be an integer. However, when $(N1)\bmod(A) > 0$, rounding that addend (or the sum of all addends) towards negative infinity would be required to avoid being negatively biased by $L$. Now, using this strategy in the equation (3) gives the following equation (4).

$$J(r',c') = floor\left(\left(\dfrac{N1+((N1)\bmod(A))}{A}\right) + \left(\dfrac{N2+((N2)\bmod(B))}{B}\right) + \left(\dfrac{N3+((N3)\bmod(C))}{C}\right) + \left(\dfrac{N4+((N4)\bmod(D))}{D}\right)\right) \quad (4)$$

where the divisors, $A = \left(\dfrac{1}{W1+L}\right), B = \left(\dfrac{1}{W2+L}\right), C = \left(\dfrac{1}{W3+L}\right)$ and $D = \left(\dfrac{1}{W4+L}\right)$. Now, let us simplify the equation (4) by generalizing it as shown in the equation (5).

$$J = floor\left(\sum_{x=1}^{n}\left(\dfrac{N_x + (N_x)\bmod(V_x)}{V_x}\right)\right) \quad (5)$$

Here, the $A, B, C, D$ divisors are represented by $V_x$ for $x \to n$. Here, $n$ can be any non-zero positive integer. $J$ represents the resulting integer. Figure 2 shows the round-off errors produced by the *floor (with floor(x) = x – xmod1)*, *ceil*, *fix*, *round* and *modulo or improved-floor* functions. The round-off errors are calculated with reference to $N_x/V_x$ where $N_x$ = [13  11  17  19  14  13  11  11  3  9] and $V_x$ = [4  10  3  8  3  5  7  9  2  6].

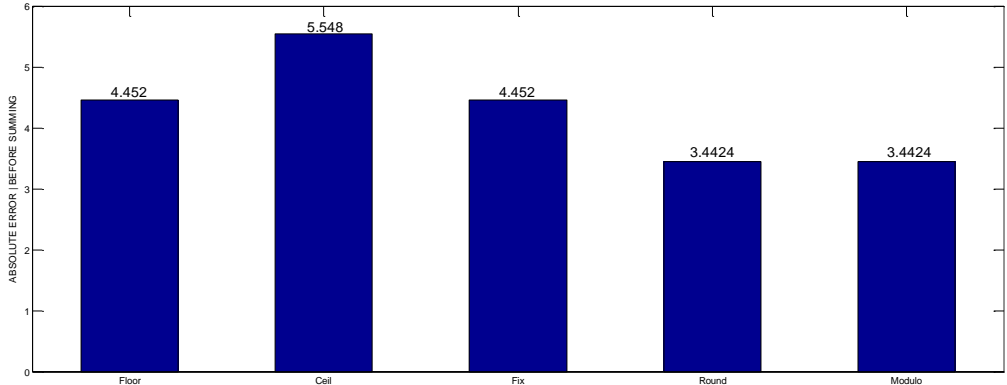

Figure 2. The round-off errors produced by each rounding scheme before adding all the addends of equation (1)

In Figure 2, the *improved-floor* produced less round-off errors than the standard *floor* function. It also performed better than the standard *ceil* and *fix* functions, except the standard *round* function where they remained tied. However, separate experiments showed that the *round* function increased the corrupting noise that negatively affected the fidelity of the images interpolated contrary to the *improved-floor* scheme. Figure 3 shows the round-off absolute errors produced by each scheme but after adding together all the addends. In this manner, only the *ceil* function produced more errors than the remaining functions. This suggests that the rounding scheme has to be appreciated according to its wide effective application rather than how it produces the least round-off errors. In conclusion, the round-off operation performed initially generates more errors than that performed later. Considering real data examples, unlike the nearest neighbour selected pixels, the bilinear interpolation pixel locations simultaneously tend to occupy more than one pixel location, which allow them to be evenly distributed over the H.R. image. To evaluate these effects, one can find the case in which pixels correlate better.

$$J(r',c') = floor\left(\begin{array}{c}\left(\dfrac{N1}{A}\right)+\left(\dfrac{N2}{B}\right)+\left(\dfrac{N3}{C}\right)+\left(\dfrac{N4}{D}\right)+\left(\dfrac{(N1)\bmod(A)}{A}\right)+\\ \left(\dfrac{(N2)\bmod(B)}{B}\right)+\left(\dfrac{(N3)\bmod(C)}{C}\right)+\left(\dfrac{(N4)\bmod(D)}{D}\right)\end{array}\right) \quad (6)$$

It is possible by exchanging only the C and D positions of the integer factor $N3$ and $N4$ of equation (6), attempting to achieve higher pixel correlation, as shown in equation (7).

$$J(r',c') = floor\left(\begin{array}{c}\left(\dfrac{N1}{A}\right)+\left(\dfrac{N2}{B}\right)+\left(\dfrac{N3}{D}\right)+\left(\dfrac{N4}{C}\right)+\left(\dfrac{(N1)\bmod(A)}{A}\right)+\\ \left(\dfrac{(N2)\bmod(B)}{B}\right)+\left(\dfrac{(N3)\bmod(C)}{C}\right)+\left(\dfrac{(N4)\bmod(D)}{D}\right)\end{array}\right) \quad (7)$$

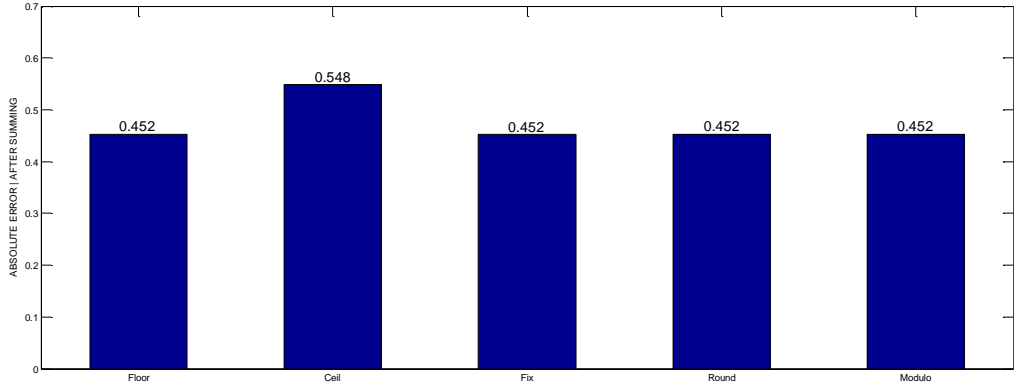

Figure 3. The round-off errors produced by each rounding scheme after adding all the addends of equation (1)

## 4. Experimental results

The exemplary standard rounding schemes have been selected to be tested against the improved-floor or modulo scheme. Equations (8) and (9) show the mathematical expressions on which the $BA_F$ and $BA_R$ have been implemented. Note that $BA_M$ is an abbreviation of the bilinear algorithm with the modulo or *improved-floor* scheme. Bilinear algorithms with the *floor* (round towards negative infinity) and *round* (rounds away from zero) are abbreviated as $BA_F$ and $BA_R$, respectively.

$$J(r',c') = floor(N1 \times W1 + N2 \times W2 + N3 \times W3 + N4 \times W4) \quad (8)$$
$$J(r',c') = round(N1 \times W1 + N2 \times W2 + N3 \times W3 + N4 \times W4) \quad (9)$$

The feature similarity index (FSIM) is one of the structural based IQA metrics used with the signal-to-noise-ratio (SNR), peak-signal-to-noise-ratio (PSNR), pixel-based IQA metrics as well as the linear correlation (CORR2).

Table 7. Top → bottom | × 2 | × 3 | × 4 | × 5 | Lena (128 ×128)

| FSIM | | | SNR | | | PSNR | | | CORR | | |
|---|---|---|---|---|---|---|---|---|---|---|---|
| $BA_F$ | $BA_R$ | $BA_M$ | $BA_F$ | $BA_R$ | $BA_M$ | $BA_F$ | $BA_R$ | $BA_M$ | $BA_F$ | $BA_R$ | $BA_M$ |
| 0.87318 | 0.87318 | 0.89739 | 19.413 | 19.413 | 19.887 | 25.164 | 25.164 | 25.637 | 0.95563 | 0.95563 | 0.96066 |
| 0.85925 | 0.85925 | 0.89726 | 17.61 | 17.61 | 18.501 | 23.355 | 23.355 | 24.246 | 0.9326 | 0.9326 | 0.94545 |
| 0.83045 | 0.83045 | 0.86152 | 16.94 | 16.94 | 17.668 | 22.621 | 22.621 | 23.349 | 0.92034 | 0.92034 | 0.9326 |
| 0.83252 | 0.83252 | 0.86134 | 16.701 | 16.701 | 17.244 | 22.443 | 22.443 | 22.986 | 0.91672 | 0.91672 | 0.92688 |

Table 7 shows the $BA_M$, $BA_F$ and $BA_R$ results at different scale ratios of Lena image as indicated in the caption. In all the cases presented, the $BA_M$ provided the results superior to those given by $BA_R$ and $BA_F$. This was due to the advantage $BA_M$ has over the $BA_R$ and $BA_F$ of adaptable divisors. Similarly, the $BA_R$ and $BA_F$ produced exactly the same results. This occurred because (as shown in Figure 3) when the addends containing non-integers were all added before being rounded both the schemes produced the same errors. Note that the accuracy in the calculations can also be negatively biased by the round-off effects of the 8-bits unsigned integer function (*UINT8*) required to display the grayscale values. In Figure 4, the (d) image interpolated using $BA_M$ appears sharper. However, it shows slight jaggedness on the edge features compared with the (b) and (d) images interpolated using the $BA_F$ (b) $BA_R$ (c). Moreover, the edge features remained equally blurred. $BA_M$'s superior results remained same also with the Baboon image (see Table 8). Furthermore, the edge features of the (d) image produced by the $BA_M$ do not appear jagged, instead tends to be blurred with slight increase in sharpness.

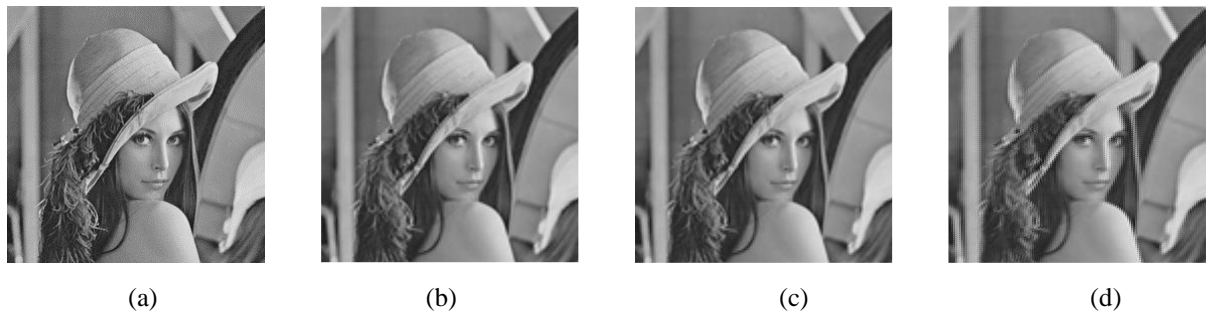

(a)  (b)  (c)  (d)

Figure 4. (a) Original image-256×256 | (b) $BA_F$ interpolated image-256×256 | (c) $BA_R$ interpolated image-256×256 | (d) $BA_M$ interpolated image-256×256

Table 8. Top → bottom | × 2 | × 3 | × 4 | × 5 | Baboon (128 ×128)

| FSIM | | | SNR | | | PSNR | | | CORR | | |
|---|---|---|---|---|---|---|---|---|---|---|---|
| $BA_F$ | $BA_R$ | $BA_M$ | $BA_F$ | $BA_R$ | $BA_M$ | $BA_F$ | $BA_R$ | $BA_M$ | $BA_F$ | $BA_R$ | $BA_M$ |
| 0.80515 | 0.80515 | 0.80874 | 17.169 | 17.169 | 17.237 | 23.319 | 23.319 | 23.386 | 0.89623 | 0.89623 | 0.89857 |
| 0.85587 | 0.85587 | 0.86294 | 16.805 | 16.805 | 16.781 | 23.089 | 23.089 | 23.064 | 0.89145 | 0.89145 | 0.89304 |
| 0.79612 | 0.79612 | 0.78953 | 16.552 | 16.552 | 16.612 | 22.77 | 22.77 | 22.83 | 0.88106 | 0.88106 | 0.88501 |
| 0.82699 | 0.82699 | 0.82903 | 16.276 | 16.276 | 16.276 | 22.559 | 22.559 | 22.559 | 0.87606 | 0.87606 | 0.87905 |

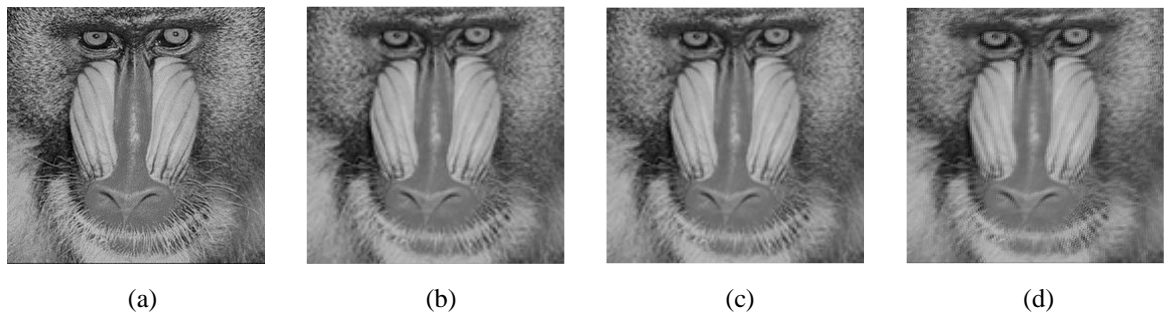

(a)  (b)  (c)  (d)

Figure 5. (a) Original image-256×256 | (b) $BA_F$ interpolated image-256×256 | (c) $BA_R$ interpolated image-256×256 | (d) $BA_M$ interpolated image-256×256

Also, as can be seen in Table 9, with the Cameraman images, the $BA_M$ output remained superior in all the scenarios presented.

Table 9. Top → bottom | × 2 | × 3 | × 4 | × 5 | Cameraman (128 ×128)

| FSIM | | | SNR | | | PSNR | | | CORR | | |
|---|---|---|---|---|---|---|---|---|---|---|---|
| $BA_F$ | $BA_R$ | $BA_M$ | $BA_F$ | $BA_R$ | $BA_M$ | $BA_F$ | $BA_R$ | $BA_M$ | $BA_F$ | $BA_R$ | $BA_M$ |
| 0.83649 | 0.83649 | 0.84962 | 17.942 | 17.942 | 18.289 | 23.622 | 23.622 | 23.968 | 0.96258 | 0.96258 | 0.96572 |
| 0.85275 | 0.85275 | 0.88203 | 16.42 | 16.42 | 17.247 | 22.101 | 22.101 | 22.928 | 0.94658 | 0.94658 | 0.95615 |
| 0.80657 | 0.80657 | 0.82554 | 15.887 | 15.887 | 16.629 | 21.504 | 21.504 | 22.245 | 0.93867 | 0.93867 | 0.94838 |
| 0.81824 | 0.81824 | 0.84131 | 15.646 | 15.646 | 16.23 | 21.32 | 21.32 | 21.904 | 0.93594 | 0.93594 | 0.94434 |

With the (d) output of $BA_M$ (see Figure 6), the image appears sharper although there seems to be a slight increase in what can be perceived as jaggedness of edge features compared to the blurred $BA_F$ and $BA_R$ cases.

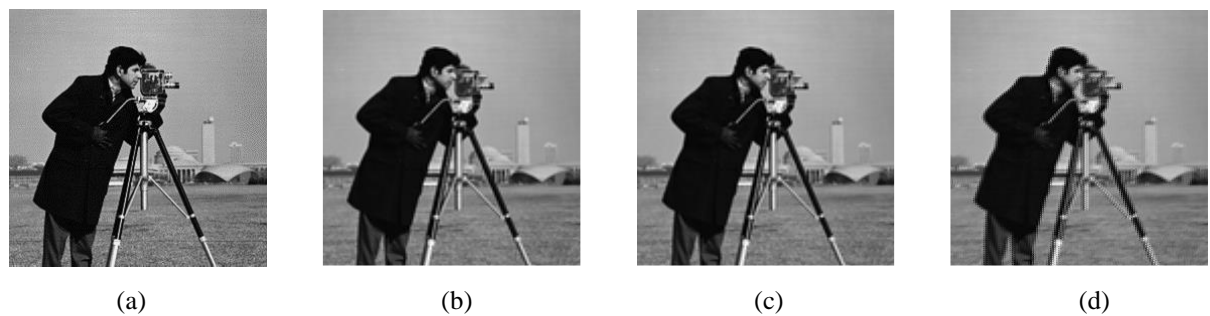

(a)       (b)       (c)       (d)

Figure 6. (a) Original image-256×256 | (b) $BA_F$ interpolated image-256×256 | (c) $BA_R$ interpolated image-256×256 | (d) $BA_M$ interpolated image-256×256

In Table 10, the effects of BAM remained positive even with the Pepper images. The IQA metrics values produced by the BAM remained higher than those of the BAF and BAR.

Table 10. Top → bottom | × 2 | × 3 | × 4 | × 5 | Peppers (128 ×128)

| FSIM | | | SNR | | | PSNR | | | CORR | | |
|---|---|---|---|---|---|---|---|---|---|---|---|
| $BA_F$ | $BA_R$ | $BA_M$ | $BA_F$ | $BA_R$ | $BA_M$ | $BA_F$ | $BA_R$ | $BA_M$ | $BA_F$ | $BA_R$ | $BA_M$ |
| 0.89 | 0.89 | 0.90785 | 19.782 | 19.782 | 20.309 | 25.793 | 25.793 | 26.32 | 0.96959 | 0.96959 | 0.9734 |
| 0.88773 | 0.88773 | 0.91846 | 17.893 | 17.893 | 18.689 | 23.9 | 23.9 | 24.696 | 0.95287 | 0.95287 | 0.96104 |
| 0.85544 | 0.85544 | 0.8804 | 17.068 | 17.068 | 17.658 | 23.011 | 23.011 | 23.6 | 0.94237 | 0.94237 | 0.94971 |
| 0.8647 | 0.8647 | 0.88675 | 16.881 | 16.881 | 17.287 | 22.885 | 22.885 | 23.291 | 0.94045 | 0.94045 | 0.94612 |

As can be seen in Figure 7, the interpolated or images produced by $BA_M$ appear only slightly sharper than by $BA_F$ and $BA_R$.

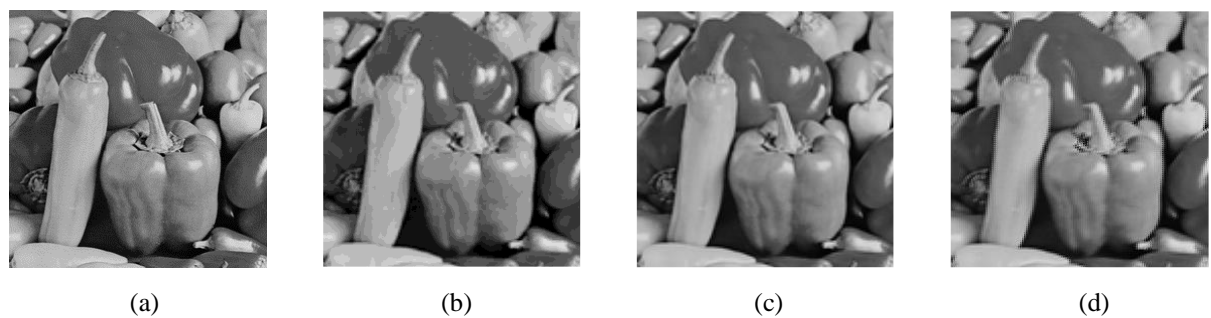

(a)       (b)       (c)       (d)

Figure 7. (a) Original image-256×256 | (b) BAF interpolated image-256×256 | (c) BAR interpolated image-256×256 | (d) BAM interpolated image-256×256

Furthermore, in Tables 11, 12 and 13, with the lake, house and girl images the BAM performances remained comparatively higher than others. In each set, the IQA metrics values obtained by the BAM remained superior to those of the BAF and BAR. The sharpness remained high, but debatable, with the images produced by BAM than those produced by both BAF and BAR (Figures 8, 9 and 10).

Table 11. Top → bottom | × 2 | × 3 | × 4 | × 5 | Lake (128 ×128)

| FSIM | | | SNR | | | PSNR | | | CORR | | |
|---|---|---|---|---|---|---|---|---|---|---|---|
| $BA_F$ | $BA_R$ | $BA_M$ | $BA_F$ | $BA_R$ | $BA_M$ | $BA_F$ | $BA_R$ | $BA_M$ | $BA_F$ | $BA_R$ | $BA_M$ |
| 0.84788 | 0.84788 | 0.86004 | 17.863 | 17.863 | 18.236 | 23.101 | 23.101 | 23.474 | 0.96183 | 0.96183 | 0.96524 |
| 0.86715 | 0.86715 | 0.88803 | 17.162 | 17.162 | 17.493 | 22.455 | 22.455 | 22.786 | 0.95563 | 0.95563 | 0.95938 |
| 0.81229 | 0.81229 | 0.82869 | 15.58 | 15.58 | 16.059 | 20.751 | 20.751 | 21.231 | 0.93458 | 0.93458 | 0.94143 |
| 0.82536 | 0.82536 | 0.84722 | 15.418 | 15.418 | 15.763 | 20.644 | 20.644 | 20.989 | 0.93277 | 0.93277 | 0.93813 |

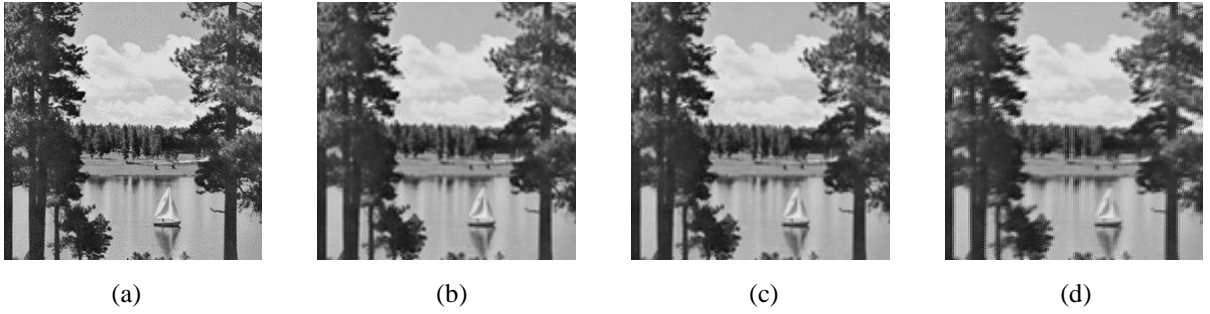

(a) (b) (c) (d)

Figure 8. (a) Original image-256×256 | (b) BF interpolated image-256×256 | (c) BAR interpolated image-256×256 | (d) BAM interpolated image-256×256

Table 12. Top → bottom | × 2 | × 3 | × 4 | × 5 | House (128 ×128)

| FSIM | | | SNR | | | PSNR | | | CORR | | |
|---|---|---|---|---|---|---|---|---|---|---|---|
| $BA_F$ | $BA_R$ | $BA_M$ | $BA_F$ | $BA_R$ | $BA_M$ | $BA_F$ | $BA_R$ | $BA_M$ | $BA_F$ | $BA_R$ | $BA_M$ |
| 0.88287 | 0.88287 | 0.90084 | 23.112 | 23.112 | 23.6 | 27.858 | 27.858 | 28.346 | 0.98415 | 0.98415 | 0.98608 |
| 0.89781 | 0.89781 | 0.92033 | 22.054 | 22.054 | 22.373 | 26.836 | 26.836 | 27.154 | 0.98026 | 0.98026 | 0.98232 |
| 0.85246 | 0.85246 | 0.87277 | 20.634 | 20.634 | 21.2 | 25.344 | 25.344 | 25.91 | 0.97127 | 0.97127 | 0.97506 |
| 0.86574 | 0.86574 | 0.88437 | 20.395 | 20.395 | 20.773 | 25.138 | 25.138 | 25.515 | 0.97014 | 0.97014 | 0.97304 |

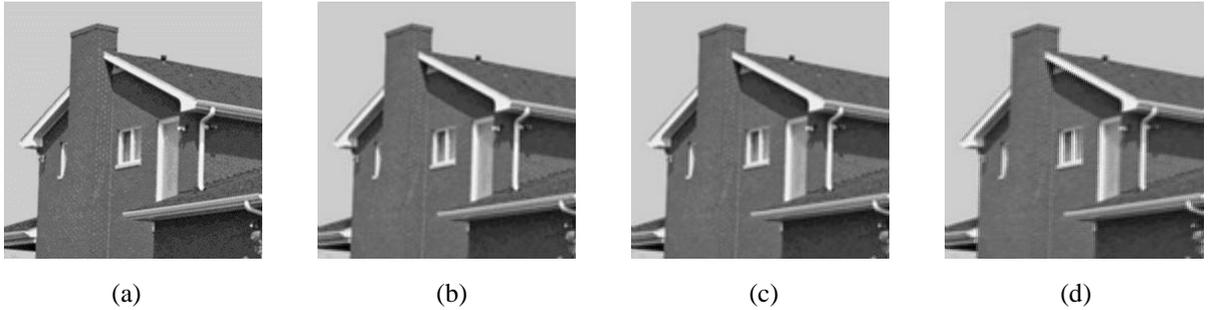

(a) (b) (c) (d)

Figure 9. (a) Original image-256×256 | (b) BAF interpolated image-256×256 | (c) BAR interpolated image-256×256 | (d) BAM interpolated image-256×256

Table 13. Top → bottom | × 2 | × 3 | × 4 | × 5 | Girl (128 ×128)

| FSIM | | | SNR | | | PSNR | | | CORR | | |
|---|---|---|---|---|---|---|---|---|---|---|---|
| $BA_F$ | $BA_R$ | $BA_M$ | $BA_F$ | $BA_R$ | $BA_M$ | $BA_F$ | $BA_R$ | $BA_M$ | $BA_F$ | $BA_R$ | $BA_M$ |
| 0.88093 | 0.88093 | 0.89681 | 20.793 | 20.793 | 21.227 | 25.981 | 25.981 | 26.415 | 0.95123 | 0.95123 | 0.95647 |
| 0.89904 | 0.89904 | 0.92046 | 20.156 | 20.156 | 20.608 | 25.404 | 25.404 | 25.856 | 0.94398 | 0.94398 | 0.95057 |
| 0.85382 | 0.85382 | 0.87288 | 17.88 | 17.88 | 18.392 | 22.991 | 22.991 | 23.503 | 0.90557 | 0.90557 | 0.91586 |
| 0.86571 | 0.86571 | 0.8843 | 17.855 | 17.855 | 18.238 | 23.027 | 23.027 | 23.409 | 0.90533 | 0.90533 | 0.9136 |

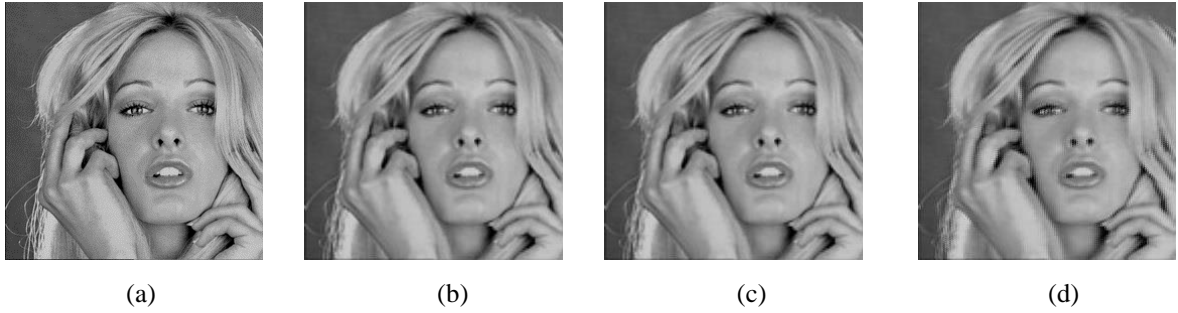

| (a) | (b) | (c) | (d) |

Figure 10. (a) Original image-256×256 | (b) BAF interpolated image-256×256 | (c) BAR interpolated image-256×256 | (d) BAM interpolated image-256×256

For short computations, the *tic–toc* commands have been used for measuring the absolute elapsed time or Matlab execution time (Martin & Bill, 2011). Table 10 shows the average run time of each scheme with the Lena, Peppers, Cameraman and Baboon images of the same resolution. Note that the run-time values can vary depending on circumstances, including the version of Matlab, computer used, and source code developed.

Table 11. Absolute elapsed time in seconds

| Image Resolution | $BA_F$ | $BA_R$ | $BA_M$ |
|---|---|---|---|
| 256 × 256 | 0.00162150 | 0.00209325 | 0.00154550 |
| 384 × 384 | 0.00148525 | 0.00294600 | 0.00147400 |
| 512 × 512 | 0.00131025 | 0.00416050 | 0.00146725 |
| 640 × 640 | 0.00136225 | 0.00538875 | 0.00144600 |

As can be seen, the average time it takes for executing the $BA_M$ algorithm is almost equal to that of $BA_F$, which is in fact the same as that of the standard bilinear interpolation algorithm. The BAR took more time than both $BA_F$ and $BA_M$ in all the scenarios.

### 4. Conclusion

The effects of *improved-floor* scheme against the standard IEEE 754–2008 rounding schemes have been studied for improvements to bilinear interpolation algorithm accuracy. The cases with relatively large or small number of half-integers in H.R image interpolation have been studied by introducing and applying the new rounding scheme based on the modulo operator. The modulo operator dealt with the whole and fractional components of the bilinear equation's addends in an effort to force them to become integer-type. Furthermore, the rounding towards negative infinity was needed due to a small positive number $L$ added to the multiplicative inverse denominator to avoid the occurrence of the indeterminate cases, while dealing with or interpolating real data. The image quality assessment metrics values produced by $BA_M$ were higher than that produced by $BA_F$ and $BA_R$ cases. Moreover, in the same scenarios presented, the $BA_M$ output images proved to look sharper instead of more blurred. Furthermore, the execution time taken by the $BA_M$ and $BA_F$ remained almost tied and shorter than that by the $BA_R$. Further efforts are to be dedicated to studying the effects of performing the rounding operations focusing separately on a pre-determined range of the pixels/values.

**Acknowledgements**

Author would like to thank the reviewers and editor for their helpful comments.